\documentclass[12pt,leqno]{article}
\usepackage{amssymb}
 \usepackage{amsthm}
 \usepackage{newtxtext}          
 \usepackage{newtxmath}
 \usepackage[T1]{fontenc}
\usepackage[utf8]{inputenc}
\usepackage[table,dvipsnames]{xcolor}
\usepackage{enumitem}
\usepackage[all]{xy}
\usepackage[pdftex,bookmarks]{hyperref}
\hypersetup{
    pdfnewwindow=true,        
    colorlinks=true,          
    linkcolor=blue,           
    citecolor=blue,           
    filecolor=magenta,        
    urlcolor=cyan,            
    pdfauthor=D\"untsch \& Dzik
}
\usepackage{geometry}
\usepackage{xspace}
\usepackage[numbers,sort&compress]{natbib}

\newcommand{\klam}[1]{\ensuremath{\langle #1 \rangle}}
\newcommand{\set}[1]{\ensuremath{\{#1\}}}
\newcommand{\z}{\emptyset}
\newcommand{\df}{\ensuremath{:=}}
\newcommand{\restrict}{\ensuremath{\mathrel{\upharpoonright}}}

\newcommand{\tand}{\text{ and }}
\newcommand{\tiff}{if and only if\xspace}
\newcommand{\aright}{``$\Rightarrow$'': \ }
\newcommand{\aleft}{``$\Leftarrow$'': \ }
\newcommand{\ua}[1]{\ensuremath{\mathop{\uparrow}#1}}
\newcommand{\Iff}{\Longleftrightarrow}
\newcommand{\Iffdf}{\overset{\mathrm{df}}{\Longleftrightarrow}}
\newcommand{\Implies}{\ensuremath{\mathrel{\Rightarrow}}}
\newcommand{\two}{\ensuremath{\mathord{\mathbf{2}}}\xspace}
\newcommand{\HA}{\ensuremath{\mathbb{HA}}\xspace}
\newcommand{\IA}{\ensuremath{\mathbb{IA}}\xspace}
\renewcommand{\O}{\ensuremath{\mathcal{O}}\xspace}
\newcommand{\B}{\ensuremath{\mathcal{B}}\xspace}
\newcommand{\F}{\ensuremath{\mathcal{F}}\xspace}
\newcommand{\Va}{\ensuremath{\mathbf{V}}\xspace}
\newcommand{\onto}{\twoheadrightarrow}
\newcommand{\into}{\hookrightarrow}
\newcommand{\wlg}{w.l.o.g.\xspace }

\DeclareMathOperator{\Eq}{\mathbf{Eq}}
\DeclareMathOperator{\Fr}{\mathtt{Fr}}
\DeclareMathOperator{\Sub}{\Lambda}

\parskip=\medskipamount
\numberwithin{equation}{section}

\theoremstyle{plain}
\newtheorem{theorem}{Theorem}[section]
\newtheorem{lemma}[theorem]{Lemma}
\newtheorem{corollary}[theorem]{Corollary}

\theoremstyle{definition}

 \title{On the preservation of unification type of Heyting algebras and interior algebras}
 \author{Ivo D\"untsch\\
Dept. of Computer Science \\
Brock University\\	
St Catharines, Ontario, 
Canada \\
\href{mailto:duentsch@brocku.ca}{duentsch@brocku.ca} \and
Wojciech Dzik \\
Institute of Mathematics \\
University of Silesia \\
Katowice, Poland \\ \href{mailto:wojciech.dzik@us.edu.pl}{wojciech.dzik@us.edu.pl}
}

\begin{document}
\thispagestyle{empty}
 \maketitle
 
 \begin{abstract}
 \noindent The purpose of this note is to shed some light on the preservation of unification type of varieties of interior algebras and varieties of Heyting algebras under the functors presented by \citet{blok76} in case the varieties involved are locally finite.
  \end{abstract}
 
 \section{Introduction}\label{sec:intro}
 
 Following the ground-breaking work of \citet{McKinsey1946}, the connections between Heyting algebras and interior algebras, respectively, closure algebras,\footnote{For the definition of Heyting algebras, interior algebras or closure algebras we refer the reader to \cite{bd74} and  \cite{blok76}.}  were investigated in detail by \citet{blok76} and \citet{bd75}. \citet{blok76} has exhibited two functors between the category $\IA$ of interior algebras and the category $\HA$ of Heyting algebras
\begin{gather*}
\O\colon \IA \to \HA, \quad \B\colon \HA \to \IA.
\end{gather*}
We will investigate the following questions:
\begin{enumerate}
\item If \Va is a variety of Heyting algebras and $L \in \Va$, what is the unification type of $L$ in \Va compared to the unification type of $\B(L)$ in $\Eq(\B[\Va])$? What is the type of \Va when compared to the type of $\Eq(\B[\Va])$?
\item If \Va is a variety of interior algebras and $A \in \Va$, what is the unification type of $A$ in \Va compared to the unification type of $\O(A)$ in $\O[\Va]$? What is the type of \Va when compared to the type of $\O[\Va]$?
\end{enumerate}
We shall answer these questions if the varieties of interior algebras involved are locally finite Grzegorczyk algebras.

 \section{Notation and first definitions}\label{sec:def}
 
 We assume that the reader is familiar with the basic concepts in lattice theory and Heyting algebras as well as interior algebras. For concepts not explained here we refer the reader to \cite{bd74} for lattice theory and Heyting algebras, \cite{blok76} for interior algebras, \cite{hs07} for category theory, and \cite{bs_ua} for universal algebra. 
  
Let $\klam{P, \leq}$ be a quasiordered set. If $x \in P$ we denote by $\ua{x} \df \set{y \in P: x \leq y}$ the order filter of $P$ generated by $x$. The classes of the equivalence relation $\theta$, defined by $x\theta y \Iffdf x \leq y \tand y \leq x$, can be partially ordered by setting $\theta(x) \leq_\theta \theta(y)$ \tiff $x \leq y$.  $M \subseteq P$ is called an \emph{antichain}, if all elements of $M$ are pairwise incomparable. $M$ is called \emph{dense}, if each element of $P$ is below some element of $M$. A \emph{$\mu$-set} is a dense antichain of $P$ \cite{ghi97}. Each $\mu$-set picks one element from each equivalence class $\theta(x)$ which is maximal with respect to $\leq_\theta$; clearly, all $\mu$-sets have the same cardinality. Following \cite{ghi97}, we say that $\klam{P,\leq}$ has unification type
\begin{quote}
\begin{description}[font=\normalfont,nosep]
\item[$1$ (\emph{unitary}),] if it has a $\mu$ set of cardinality $1$,
\item[$\omega$ (\emph{finitary}),] if it has a finite $\mu$-set of cardinality greater than $1$,
\item[$\infty$ (\emph{infinitary}),] if it has an infinite $\mu$-set,
\item[$0$ (\emph{nullary}),] if no $\mu$-set exists.
\end{description}
\end{quote}
If $L$ is a bounded distributive lattice, we let $\Fr(L)$ be its free Boolean extension;\footnote{Also called minimal Boolean extension \cite{ner59}.} we suppose \wlg that $L$ is a bounded sublattice of $\Fr(L)$. $\Fr(L)$ may be taken as the Boolean subalgebra of the power set of prime ideals of $L$, generated by the Stone embedding of $L$, see \cite[V.4, Theorem 2]{bd74}. $\Fr(L)$ is unique in the following sense: If $L$ is a bounded sublattice of a Boolean algebra $B$, then the Boolean subalgebra of $B$ generated by $L$ is isomorphic to $\Fr(L)$; in particular, $B \cong \Fr(L)$ if $L$ generates $B$.

If $K$ is a class of algebras of the same type, $\Eq(K)$ denotes the variety generated by $K$. We will only consider non-trivial varieties. For a variety \Va we denote by $\Sub(\Va)$ the lattice of its subvarieties. Varieties will also be considered as categories with the usual terminology. A functor $F\colon \Va \to \Va'$ is called \emph{faithful}, if it is injective on sets of homomorphisms, and it is called \emph{full}, if it is surjective on sets of homomorphisms. $\Va$ and $\Va'$  are called \emph{categorically equivalent}, if there is a full and faithful functor $F\colon \Va \to \Va'$, and for every $B \in \Va'$ there is some $A \in \Va$ such that $F(A) \cong B$, see \cite{alb96}. If $\Va'$ is a subvariety of $\Va$, then $\Va'$ is a full subcategory of \Va, i.e. $\hom_\Va(A,B) = \hom_{\Va'}(A,B)$ for all $A,B \in \Va'$.

The varieties of Heyting algebras and interior algebras are denoted by, respectively, $\HA$ and $\IA$. The following class of interior algebras plays a decisive role in the study of the connection between Heyting algebras and interior algebras: $\klam{B,g} \in \IA$ is a \emph{Grzegorczyk algebra}, if it satisfies
\begin{gather}\label{Grz0}\tag{\textbf{Grz}}
g(x+g(x \cdot -g(x))) \leq x.
\end{gather}

 \section{Algebraic unification}
 
Algebraic unification was studied, among others, by \citet{ghi97}, and we outline his approach as we need it. Background and further details can be found in \cite{ghi97,gs04} and the references therein.\footnote{See also the related universal algebraic approach outlined by \citet{bur01}.}

Suppose that $\Va$ is a variety, and that $A$ is a finitely presented algebra relative to \Va. A \emph{unifier for $A$ in \Va} is a pair $\klam{u,B}$ where $B \in \Va$ is finitely presented and projective in $\Va$, and $u\colon A \to B$ is a homomorphism. Note that unifiability and unification are defined for finitely presented algebras only, and thus they depend on $\Va$. Finite presentability and projectivity are dependent on the class considered. 

The collection of unifiers of $A$ in $\Va$ is denoted by $U_A^\Va$, and we call $A$ \emph{unifiable relative to \Va}, if $U_A^\Va \neq \z$. If $\Va$ is understood we shall usually omit the superscript and call $A$ simply \emph{unifiable}. 

The following lemma and its corollary give a simple criterion for unifiability:
\begin{lemma}\label{lem:min}
Suppose that $M$ is a minimal algebra in the variety \Va and projective in \Va. Then, $A \in \Va$ is unifiable \tiff $M$ is a homomorphic image of $A$.
\end{lemma}
\begin{proof}
\aright Suppose that $B$ is projective in \Va, and $f:A \to B$ a homomorphism. Since $B$ is projective in \Va, $M$ is a homomorphic image of $A$ by \cite[Corollary 3.11]{Citkin2018}. Furthermore, $f[A]$ is a subalgebra of $B$, and therefore, $M$ is a homomorphic image of $f[A]$ by \cite[Proposition 3.10]{Citkin2018}. It follows that $M$ is a homomorphic image of $A$.

\aleft This follows from the fact that $M$ is a homomorphic image of $A$. 
\end{proof}
\begin{corollary}\label{cor:min}
If \Va is a variety with a bounded lattice reduct, then $A \in \Va$ is unifiable \tiff \two is a homomorphic image of $A$.
\end{corollary}
\begin{proof}
Just observe that \two is a projective minimal algebra in \Va.
\end{proof}
This shows immediately that a simple algebra with more than two elements is not unifiable, in particular, any discriminator algebra with at least four elements, since discriminator algebras are simple.

\begin{figure}[htb]
\caption{Quasiordering unifiers}\label{fig:orduni}
$$
\xymatrix{
& {A} \ar[ld]_{u} \ar[rd]^{v} \\
{B} \ar@{-->}[rr]^{h}_{\succcurlyeq} && {C} }
$$
\end{figure}
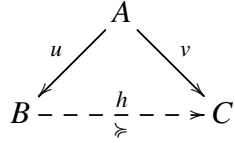

In the rest of this section we fix a variety $\Va$ with respect to which unifiers are taken. Given two unifiers $\klam{u,B}$ and $\klam{v,C}$ for $A$, we say that \emph{$\klam{u,B}$ is more general  than $\klam{v,C}$}, written as $\klam{u,B}\succcurlyeq \klam{v,C}$, if there is a homomorphism $h\colon B \to C$ such that the diagram in Figure \ref{fig:orduni} commutes, i.e. that $v = h \circ u$.

Clearly, $\succcurlyeq$ is a quasiorder. If $A$ itself is projective, then $\klam{i,A}$  is more general than any of its unifiers. If $\klam{u,B} \preccurlyeq \klam{v,C}$ and  $\klam{u,B} \succcurlyeq \klam{v,C}$ we write $\klam{u,B} \approx \klam{v,C}$. The relation $\approx$ is an equivalence relation on $U_A$, and $U_A/\approx$ can be partially ordered as described above. We define the \emph{unification type $t^\Va(A)$ of $A$ in \Va} as the type of the quasiordered set  $\klam{U_A, \succcurlyeq}$ as defined in Section \ref{sec:def}. If \Va is understood, we shall omit the superscript.

Generalizing unifier type to varieties, we say that \Va has unification type
\begin{description}[font=\normalfont,nosep]
\item[\emph{unitary}] if every unifiable $A \in \Va$ has type $1$.
\item[\emph{finitary}] if every unifiable $A \in \Va$ has type $1$ or type $\omega$, and there is some unifiable $A \in \Va$ with type $\omega$,
\item[\emph{infinitary}] if every unifiable $A \in \Va$ has type $1, \omega$ or $\infty$, and there is some unifiable $A \in \Va$ with type $\infty$,
\item[\emph{nullary}] if there is some $A \in \Va$ with type $0$.
\end{description}
The unification type of $\Va$ is denoted by $t(\Va)$. If we order unification types by $1 \lneq \omega \lneq \infty \lneq 0$, we see that $t(\Va) = \max\set{t(A): A \in \Va}$, see e.g. \cite[Definition 3.4]{bs01}.

\section[Interior algebras and Heyting algebras]{Interior algebras and Heyting algebras}\label{sec:HA}
An important role in the study of the connection between $\IA$ and $\HA$ are played by $^*$-algebras. Suppose that $\Va$ is a variety of interior algebras. Following \citet[Definition I.2.14]{blok76} we call $B \in \Va$ a \emph{$^*$-algebra}, if it is generated by its open elements. We denote by $\Va^*$ the variety generated by the $^*$-algebras it contains. If $\Va = \Va^*$, we call it a \emph{$^*$-variety}. The variety $\IA^*$ can be neatly described:
\begin{theorem}\label{thm:Grz*}\cite[Corollary III.7.9]{blok76}
$\IA^*$ is the variety of Grzegorczyk algebras.
\end{theorem}
This implies that each subvariety of $\IA^*$ is also a $^*$-variety. All finite algebras in $\Va^*$ are $^*$-algebras, but not every algebra in $\Va^*$ is necessarily a $^*$-algebra \cite[Section II.2.9]{blok76}. However,
\begin{theorem}\label{thm:1.8}
\cite[Corollary III.4.5]{blok76}\label{lem loc fin} All algebras in $\Va^*$ are $^*$-algebras \tiff $\Va^*$ is locally finite.    
\end{theorem}
This can be slightly rephrased as follows:
\begin{corollary}\label{cor:lf*}
Every $A \in \Va$ is a $^*$ algebra \tiff \Va is a locally finite $^*$-variety. 
\end{corollary}
\begin{proof}
\aright If every $A \in \Va$ is a $^*$-algebra, then \Va is a $^*$-variety and \Va is locally finite by the ``$\Rightarrow$'' direction of Theorem \ref{thm:1.8}. . 

\aleft Since \Va is a $^*$-variety we have $\Va = \Va^*$, and the rest follows directly from the ``$\Leftarrow$'' direction of Theorem \ref{thm:1.8}.
\end{proof}

The next lemma shows that $A^*$ is never a retract of $A$ if it is a proper subalgebra of $A$:
\begin{lemma}\label{lem:retract}
Let $\klam{A,g} \in \IA$.  If $A^*$ is a retract of $A$, then $A = A^*$.
\end{lemma}
\begin{proof}
Suppose that $p\colon A \onto A^*$ is a retraction, and let $F \df p^{-1}(1)$. Then, 
\begin{xalignat*}{2}
a\in F &\Implies g(a) \in F, && \text{since $F$ is an open filter,} \\
 &\Implies p(g(a)) = 1, && \text{by definition of $F$,}\\
&\Implies g(a) = 1, &&\text{since $g(a) \in A^*$ and $p\restrict A^*$ is the identity}, \\
&\Implies a = 1, &&\text{since $g(a) \leq a$},\\
&\Implies F = \set{1}, \\
&\Implies p \text{ is injective.}
\end{xalignat*}
Since $p$ is also onto, it follows that $p$ is an isomorphism, hence, $A = A^*$.
\end{proof}
Following \citet{blok76}, we will now introduce the functors between \IA and \HA. For $\klam{B,g} \in \IA$ we let $B^\circ$ be its Heyting algebra of open elements, that is, the $0,1$-sublattice of $B$ generated by the open elements, and implication defined by $a \implies b\df g(-a + b)$ \cite[Theorem I.2.3]{blok76}. If $K \subseteq \IA$, we set 
\begin{gather}\label{def:Kcirc}
K^\circ \df \set{B^\circ: B \in K} \subseteq \HA.
\end{gather}
If $A,B \in \IA$ and $h\colon A \to B$ is a homomorphism, we define $h^\circ \df h\restrict A^\circ$. It is well known that $h^\circ$ is an \HA-homomorphism, and that $h^\circ$ is onto if $h$ is onto \cite[Theorem I.2.11]{blok76}.

Next, we define the functor $\O\colon \IA \to \HA$. Let $A,B \in \IA$ and $h\colon A \to B$ be a homomorphism. 
\begin{align}
\O(A) &\df A^\circ. \\
\O(h) &\df h^\circ.
\end{align}
\begin{theorem}\label{thm:O}
\begin{enumerate}
\item \cite[Corollary I.2.12]{blok76} $\mathcal{O}$ is a covariant functor between the categories of interior algebras and Heyting algebras, and the diagram of Figure \ref{fig:O} commutes.
\end{enumerate}
\end{theorem}
\begin{figure}[h!tb]
\caption{The functor $\O$}\label{fig:O}
$$
\xymatrix{
A \ar@{->}[d]_f \ar@{->}[r]^\O & A^\circ \ar@{->}[d]^{\O(h) = f^\circ}\\
B \ar@{->}[r]_\O & B^\circ
}
$$
\end{figure}

We extend $\O$ over varieties in the following way: Let $\gamma\colon \Sub(\IA) \to \Sub(\HA)$ be defined by
\begin{gather}\label{def:gamma}
\gamma(\Va) \df \Va^\circ.
\end{gather}
By the definitions, $\gamma(\Va) = \set{A^\circ: A \in \Va} = \set{\O(A): A \in \Va} = \O[\Va]$. Furthermore, $\gamma(\Va) = \gamma(\Va^*)$ since $\O(A) = O(A^*)$ for each $A \in \Va$.
\begin{theorem}\label{thm:gamma}\cite[Theorem III.7.3]{blok76} 
$\gamma$ is a complete surjective lattice homomorphism.
\end{theorem}
Conversely we define the functor \B from \HA to \IA. For $L \in \HA$ its free Boolean extension $\Fr(L)$ can be augmented by an operator $g_L$ defined as follows: If $a \in \Fr(L)$, $a = (-u_i + v_1) \cdot \ldots \cdot (-u_n+v_n)$ for some $u_i, v_i \in L$, we set
\begin{gather}\label{def:g}
g_L(a) \df (u_i \implies v_1) \cdot \ldots \cdot (u_n \implies v_n),
\end{gather}
\begin{lemma}\label{lem:B} 
Let $L \in \HA$ and $B \in \IA$. 
\begin{enumerate}
\item \cite[Theorem I.2.13]{blok76} $g_L$ as defined by \eqref{def:g} is the unique interior operator on $\Fr(L)$ such that $\Fr(L)^\circ \cong L$. 
\item \cite[Theorem I.2.15]{blok76} Let $h:L \to B^\circ$ be an \HA-homomorphism. Then, there is a unique \IA-homomorphism $\overline{h}\colon \klam{\Fr(L),g_L} \to B$ such that $\overline{h}^\circ = h$.
\end{enumerate}
\end{lemma}
Now we can define the functor $\B\colon \HA \to \IA$. If $L, M\in \HA$, and  $h\colon L \to M$ is a homomorphism, we set
\begin{align}
&\B(L) \df  \klam{\Fr(L),g_L}, \\
&\B(h) \df \overline{h}. \label{Bh}
\end{align}
\begin{theorem}\label{thm:BV} \cite[Corollary I.2.17]{blok76} The assignment $\B\colon \HA \to \IA$ is a covariant functor which preserves  injective and surjective morphisms, and the diagram of Figure \ref{fig:B} commutes. Furthermore, $\B$ is a full embedding.
\end{theorem}

\begin{figure}[h!tb]
\caption{The functor $\B$}\label{fig:B}
 $$
\xymatrix{
L \ar@{->}[d]_h \ar@{->}[r]^\B & \Fr(L) \ar@{->}[d]^{\B(h) = \overline{h}} \\
L_1 \ar@{->}[r]_\B & \Fr(L_1)
}
$$
\end{figure}

\begin{lemma}\label{lem:CatEquiv}
If $\Va \leq \HA$, then \Va and $\B[\Va]$ are categorically equivalent.
\end{lemma}
\begin{proof} Since $\B\restrict\Va\colon \Va \to \B[\Va]$ is bijective, it remains to show that $\B\restrict\Va$ is fully faithful. This follows from the facts that  $\B$ is fully faithful by Theorem \ref{thm:BV}, and that $\Va$ is a full subcategory of \HA.
\end{proof}

The next result establishes a connection between $\B,\O$ and $^*$-algebras:
\begin{theorem}\label{thm:BGrz}
\begin{enumerate}
\item Suppose that $A,B \in \IA$ and $h\colon A \to B$ is a homomorphism. Then, $B^* \cong \B(\O(B))$ and $h\restrict B^* = \B(\O(h))$.
\item Suppose that $\Va$ is a variety of Heyting algebras. Then, $\Eq(\B[\Va]) = \Eq(\B[\Va])^*$.
\item Suppose that \Va is a variety of interior algebras. Then, $\Eq(\B[\O[\Va]]) = \Va^*$.
\end{enumerate}
\end{theorem}
\begin{proof}
1. $\B(\O(B)) = \B(B^\circ)$ is the free Boolean extension of $\B^\circ$ which is isomorphic to $B^{*}$. Furthermore, $\B(\O(h)) = \B(h^\circ) = \overline{h^\circ} = h^*$.

2. $\Eq(\B[\Va])^*$ is the variety generated by the $^*$-algebras contained in $\Eq(\B[\Va])$. Since $\Eq(\B[\Va])$ is generated by a collection of $^*$-algebras, the claim follows.

3. $\Eq(\B[\O[\Va]])$ is the equational class generated by the $^*$-algebras in \Va, hence, it is equal to $\Va^*$.
\end{proof}

In analogy to $\gamma$ we define a mapping $\rho\colon \Sub(\HA) \to \Sub(\IA)$ by 
\begin{gather}\label{def:rho}
\rho(\Va) \df \set{A \in \IA: \O(L) \in \Va}.
\end{gather}. 
\begin{theorem}\label{thm:rho}\cite[Theorem III.7.4]{blok76} 
$\rho$ is an injective homomorphism of bounded lattices, and $\gamma \circ \rho$ is the identity on $\Sub(\HA)$. 
\end{theorem}

$\rho(\Va)$ is the largest variety $\mathbf{W}$ of closure algebras such that $\O(A) \in \Va$ for all $A \in \mathbf{W}$ \cite[p.90]{blok76}, in other words, that $\mathbf{W}^\circ = \Va$. Its corresponding logic is the smallest modal companion of the logic induced by $\Va$. Note that Theorem \ref{thm:rho}  implies that $\O(\B(L)) \cong L$ for all $L \in \HA$. 

Next we turn to the smallest variety $\mathbf{W} \leq \IA$ with the property that $\mathbf{W}^\circ = \Va$.  Define $\rho^*\colon \Sub(\HA) \to \Sub(\IA^*)$ by
\begin{gather}\label{def:rho*}
\rho^*(\Va) \df \Eq(\B[\Va]).
\end{gather}
Note that $(\gamma \circ \rho^*)(\Va) = (\gamma \circ \rho)(\Va)$, since $\rho(\Va)^* = \rho^*(\Va)$.

The next result implies the famous Blok-Esakia Theorem:
\begin{theorem}\label{thm:rho*}\cite[Theorem III.7.10]{blok76}
$\rho^*$ is a lattice isomorphism.
\end{theorem}
In the rest of this section we establish several results relating projectivity of algebras in $\Va \leq \IA$ to projectivity in $\Va^*$ that will be used in later sections. A central result is the following observation by \citeauthor{blok76}:
\begin{lemma}\label{thm:7.16} \cite[Theorem I.7.16]{blok76} 
If $\Va$ is a variety of interior algebras, then a $^*$-algebra $B \in \Va$ is projective in $\Va$ \tiff $\O(B)$ is projective in $\gamma(\Va)$.
\end{lemma}

\begin{lemma}\label{lem:epi}
Suppose that  $A,B \in \Va\leq\IA$ and $h\colon A \to B$ is a homomorphism. Then, $h^* \df h\restrict A^*\colon A^* \to B^*$ is a homomorphism, $h^*$ is onto, if $h$ is onto, and $h^*$ is one-one, if $h$ is one-one.
\end{lemma}
\begin{proof}
The functor $\B\circ \O: \Va \to \Va^*$ defined by $(\B \circ\O)(A) = \Fr(A^\circ)) \cong A^*$ and $(\B \circ \O)(h) = \overline{(h^\circ)} = h^*$ preserves surjectivity and injectivity of morphisms, see \cite[Corollaries I.2.12, I.2.17]{blok76}
\end{proof}
This lemma has several useful consequences:
\begin{lemma}\label{lem:retract2}
If $A,B \in \Va$ and $p\colon A \onto B$ is a retraction witnessed by $q\colon B \into A$, then $p^*$ is a retraction witnessed by $q^*$.
\end{lemma}
\begin{proof}
This follows from Lemma \ref{lem:epi} and the fact that the functors preserve composition of morphisms. 
\end{proof}
Lemma \ref{lem:retract2} shows that if $B$ is a retract of $A$ in \Va, then $B^*$ is a retract of $A^*$ in $\Va^*$. We don't know whether in general projectivity of $B$ in \Va implies projectivity of $B^*$ in $\Va^*$. The situation is pictured in the following diagram:
$$
\xymatrix
{  
B  \ar@/^/[d]^p & B^* \ar@/^/[d]^{p^*} \ar@{>->}[l] \\
A  \ar@/^/[u]^q & A^* \ar@/^/[u]^{q^*} \ar@{>->}[l] 
}
$$
Conversely, if $A^*$ is a retract of $B^*$, then $A$ is not necessarily a retract of $B$: Suppose that  \Va is the variety of monadic algebras. If $A,B$ are different discriminator algebras with at least four elements, then $A^* = B^* = \two$, and $A$ is not a retract of $B$. Furthermore, $B^* = \two$ is projective in $\Va^*$, hence in \Va, but $B$ is not projective in \Va, see \cite[Corollary 6.2]{kq76}.

We have the following weaker result:

\begin{lemma}\label{lem:proj*}
If $B^*$ is projective in $\Va^*$ \tiff $B^*$ is projective in \Va.
\end{lemma}
\begin{proof}
\aright Suppose that $A,B^* \in \Va$ and that $f\colon A \onto B^*$ is an epimorphism. Then, $f^*\colon A^* \onto B^*$ is an epimorphism, and projectivity of $B^*$ in $\Va^*$ implies that there is some $g\colon B^* \into A^*$ such that $f^* \circ g$ is the identity on $B^*$. Since $f^* \circ g = f \circ g$ the claim follows.

\aleft If $B^*$ is projective in \Va, then it is projective in $\Va^*$ since projectivity is hereditary with respect to subvarieties.
\end{proof}

\section{Preservation of unification type}\label{sec:preserv}

In this section we will explore the preservation of unification type by the functors \O and \B of algebras, respectively, by the operators $\gamma, \rho$ and $\rho^*$ of varieties. With some abuse of notation we will regard the functors $\O$ and $\B$ as mappings from unifiers as well, so that $\B(u,L)$ means the pair $\klam{\B(u),\B(L)}$, and analogously for $\O$.  We first recall a result by \citeauthor{alb96}, which is the only previous result we know on preservation of unification type:
\begin{theorem}\label{lem:Aequiv}\cite[Theorem 1]{alb96} 
Suppose that the varieties $\Va, \Va'$ are categorically equivalent via a functor \F.  Then, the unification type of each algebra is preserved by \F. Consequently, $t(\Va) = t(\Va')$.
\end{theorem}

To show that the image of a unifier under a functor is the a unifier in the image, we need to show that finite presentability and projectivity are preserved in the respective varieties. Finite presentation is taken care of by the following result:

\begin{theorem}\label{thm:fp}\cite[Theorem 4.8]{cit12}
\begin{enumerate}
\item If $\Va$ is a variety of interior algebras and $A \in \Va$ is finitely presented in \Va, then $\O(A)$ is finitely presented in $\O[\Va]$.
\item If $\Va$ is a variety of Heyting algebras and $L \in \Va$ is finitely presented in \Va, then $\B(L)$ is finitely presented in $\rho^*(\Va)$.
\end{enumerate}
\end{theorem}

\subsection{From \HA to \IA}\label{sec:HAtoIA}

In this section we suppose that \Va is a variety of Heyting algebras. Our first aim is to show that $\B$ preserves unifiers:
\begin{lemma}\label{lem:Vproj}
If $L \in \Va$ is projective in \Va, then $\B(L)$ is projective in $\rho(\Va)$. 
\end{lemma}
\begin{proof}
By Theorem \ref{thm:BGrz}  $\B(L)$ is a $^*$-algebra, and now
\begin{gather*}
\B(L) \text{ is projective in }\rho(\Va) \Iff \O(\B(L)) \text{ is projective in } \gamma(\rho(\Va))
\end{gather*}
by Lemma \ref{thm:7.16}. The claim now follows from the fact that $\O(\B(L)) = L$, and $L$ is projective in $\gamma(\rho(\Va)) = \Va$.
\end{proof}
Since projectivity is inherited by subvarieties and $\B(L) \in \Va^*$, projectivity of $L$ in \Va implies projectivity of $\B(L)$ in $\rho^*(\Va)$. Theorem \ref{thm:fp} and Lemma \ref{lem:Vproj} immediately imply preservation of unifiers:
\begin{theorem}\label{thm:presproj}
If $\klam{u,M}$ is a unifier of $L$ in $\Va$, then $(\klam{\B(u), \B(M)}$ is a unifier of $\B(L)$ in $\rho^*(\Va)$.
\qed\end{theorem}

If $L \in \Va$ is unifiable, we let 
\begin{gather}\label{def:tau}
\tau\colon U^\Va_L \to U_{\B(L)}^{\rho^*(\Va)}, \ \tau(\klam{u,M}) \df \klam{\B(u), \B(M)}.
\end{gather}
$\tau$ is well defined by Theorem \ref{thm:presproj}. 

We now turn to the preservation of the quasiorder of the unifier set by $\tau$. The following observations is crucial:

\begin{lemma}\label{lem:locfin1}
$\rho^*(\Va)$ is locally finite \tiff $\rho^*(\Va) = \B[\Va]$.
\end{lemma}
\begin{proof}
\aright Let $B \in \rho^*(\Va))$. Since $\rho^*(\Va))$ is a locally finite $^*$ variety, $B = B^*$ by Theorem \ref{thm:1.8}. Furthermore, Theorem \ref{thm:BGrz}(1) implies that $B^* = \B(\O(B))$. It follows that $B \in \B[\Va]$. The other inclusion is obvious.

\aleft If  $\rho^*(\Va) = \B[\Va]$, then every $B \in \rho^*(\Va)$ is a $^*$-algebra and therefore, the $^*$-variety $\rho^*(\Va)$ is locally finite by Theorem \ref{thm:1.8}.
\end{proof}
\begin{theorem}\label{thm:pres}
Suppose that $L \in \Va$ is unifiable; then, $\tau$ is injective and preserves $\succcurlyeq$. Furthermore, if $\rho^*(\Va)$ is locally finite, then $\tau$ is surjective.
\end{theorem}
\begin{proof}
The mapping $\tau$ is injective, since the free Boolean extension of $L$ is unique. Let $\klam{u_0,L_0}$ and $\klam{u_1,L_1}$ be unifiers of $L$ in \Va, and $\klam{u_0,L_0} \succcurlyeq \klam{u_1,L_1}$, pictured in the commuting diagram below:
$$
\xymatrix{
& {L} \ar[ld]_{u_0} \ar[rd]^{u_1} \\
{L_0} \ar@{-->}[rr]^{v}_{\succcurlyeq} && {L_1} }
$$
Since $\B$ is a functor, the diagram
$$
\xymatrix{
& {\B(L)} \ar[ld]_{\B(u_0)} \ar[rd]^{\B(u_1)} \\
{\B(L_0)} \ar@{-->}[rr]^{\B(v)}_{\succcurlyeq} && {\B(L_1)} }
$$
commutes, and shows that $\klam{\B(u_0),\B(L_0)} \succcurlyeq \klam{\B(u_1), \B(L_1)}$. If $\rho^*(\Va)$ is locally finite, then, $\rho^*(\Va) = \B[\Va]$ by Lemma \ref{lem:locfin1}. Since $\B$ is a full functor, $\tau$ is onto.
\end{proof}

\begin{theorem}\label{thm:preslocfin1}
If $\Va \leq \HA$ and $\Eq(\B[\Va])$ is locally finite, then $t^{\Va}(L) = t^{\rho^*(\Va)}(\B(L))$ for all $L \in \Va$.
\end{theorem}
\begin{proof}
This follows from the fact that $\tau$ is bijective and preserves $\succcurlyeq$.
\end{proof}
For the preservation of the type by the functor $\rho^*$ we can utilize Theorem \ref{lem:Aequiv}:
 
 \begin{theorem}\label{thm:preslocfin}
If $\Va \leq \HA$ and $\rho^*(\Va)$ is locally finite, then $t(\Va) = t(\rho^*(\Va))$. \qed
\end{theorem} 
\begin{proof}
If $\rho^*(\Va)$ is locally finite, then $\rho^*(\Va) = \B[\Va]$ by Lemma \ref{lem:locfin1}. Since $\Va$ and $\B[\Va]$ are categorically equivalent by Lemma \ref{lem:CatEquiv}, we have $t(\Va) = t(\rho^*(\Va))$ by Theorem \ref{lem:Aequiv}.
\end{proof}

\subsection{From \IA to \HA}\label{sec:IAtoHA}
Now we turn to the direction from varieties of interior algebras to varieties of Heyting algebras. In the sequel suppose that \Va is variety of interior algebras. Let $A \in \Va$ be unifiable, and $\klam{u,B} \in U_A^\Va$. Then, $B$ is finitely presentable and projective in $\Va$. By Lemma \ref{thm:fp}, $\O(B)$ is finitely presentable in $\gamma(\Va)$, and we need to establish whether $\O(B)$ is projective in $\gamma(\Va)$. There is a close connection between projectivity of $\O(B)$ in $\gamma(\Va)$ and projectivity of $B^*$ in $\Va^*$:
\begin{lemma}\label{lem:projequiv}
Suppose that $B \in \Va$. Then,  $\O(B)$ is projective in $\gamma(\Va)$ \tiff $B^*$ is projective in $\Va^*$.
\end{lemma}
\begin{proof}
\aright Let $\O(B)$ be projective in $\gamma(\Va)$. Then, $\B(\O(B))$ is projective in $\rho(\gamma(\Va))$ by Lemma \ref{lem:Vproj}. The claim now follows from the facts that $\B(\O(B)) \cong B^*$ and $\rho^*(\gamma(\Va)) = \Va^*$.

\aleft Suppose that $B^*$ is projective in $\Va^*$. Then,  $B^*$ is projective in \Va by Lemma \ref{lem:proj*}, and therefore, $\O(B^*)$ is projective in $\gamma(\Va)$ by Lemma \ref{thm:7.16}. Since $\O(B) = \O(B^*)$, the claim follows.
\end{proof}
\begin{lemma}\label{lem:extend*}
Suppose that $A,B \in \IA$, then $\O[\hom(A,B)] = \hom(\O(A), \O(B))$ \tiff every \IA-homomorphism $h\colon A^* \to B^*$ has an extension $A \to B$.
\end{lemma}
\begin{proof}
\aright Suppose that  $h\colon A^* \to B^*$ is an \IA-homomorphism; then, $v\colon (A^*)^\circ \to (B^*)^\circ$ with $v \df h^\circ$ is an \HA-homomorphism, and $h = \overline{v}$, where $\overline{v}\colon A^* \to B^*$ is the unique homomorphism such that $\overline{v}^\circ = v$. Since $\O(A^*) = \O(A)$ and $\O(B^*) = \O(B)$ and by the hypothesis, there is some $u\colon A \to B$ such that $\O(u) = u^\circ = v$. Now, $u^* = \B(\O(u)) = \B(u^\circ) = \overline{v} = h$, and $u$ is the desired extension.

\aleft  Let $v\colon \O(A) = A^\circ \to B^\circ = \O(B)$ be an \HA - homomorphism. Then, $\overline{v}\colon A^* \to B^*$ has an extension $u\colon A \to B$ such that $u^* = \overline{v}$. Now, $\O(u) = u^\circ = (u^*)^\circ = v$.
\end{proof}
Lemma \ref{lem:extend*} can be viewed as a sort of local fullness. Next we shall establish a condition for which varieties \Va the restriction of  \O to \Va is a full functor:

\begin{theorem}\label{thm:full}
Suppose that \Va is a variety of interior algebras. Then, $\O\restrict \Va$ is full \tiff \Va is a locally finite $^*$-variety. In this case, $\Va$ and $\gamma(\Va)$ are categorically equivalent.
\end{theorem}
\begin{proof} 
\aright Let $A \in \Va$. By Corollary \ref{cor:lf*} it is enough to show that $A$ is a $^*$ algebra. Assume that $A^* \lneq A$. Let $e\colon \O(A) \to \O(A^*) = \O(A)$ be the identity; we will show that $e \not\in  O[\hom(A,A^*)]$. Assume that $h\colon A \to A^* \in \hom(A,A^*)$ such that $\O(h) = e$. Since $\B(\O(A)) = A^*$ and $\B(\O(h)) = h \restrict A^*$ is the identity on $A^*$, $h$ is a retraction, contradicting Lemma \ref{lem:retract} and $A^* \neq A$.

\aleft Suppose that $A,B \in \Va$, and let $v\colon \O(A) = A^\circ \to B^\circ = \O(B)$ be an \HA - homomorphism. Then, $\overline{v}\colon A = A^* \to  B^* = B$ is an \IA-morphism which satisfies $\O(\overline{v}^\circ) = v$  by Lemma \ref{thm:fp} and  Lemma \ref{lem:B}.

For the rest of the proof we suppose that \Va is a locally finite $^*$-variety. First, observe that $\gamma(\Va) = \O[\Va]$. Since $B = B^*$ for all $B \in \Va$ by Corollary \ref{cor:lf*}, $\O\restrict \Va$ is injective. Noting that $\O$ is surjective, we obtain that $\O\restrict \Va$ is bijective on objects. Since $\O\restrict \Va$ is full by the hypothesis, all that remains to show is that $\O\restrict \Va$ is faithful. Suppose that $A,B \in \Va$, and that $f,g\colon A \to B$ are homomorphisms. Let $\O(f) = \O(g)$.  Now, $\O(f) = f^\circ = f\restrict \Va^*$, similarly, $\O(g) = f \restrict \Va$. Since every algebra of $\Va$ is a $^*$-algebra, the claim follows.
\end{proof}
\begin{lemma}\label{lem:projrestrict}
Suppose that $\O\restrict \Va$ is a full functor, and that $B$ is projective in \Va. Then, $\O(B)$ is projective in $\gamma(\Va)$.
\end{lemma}
\begin{proof}
Using Lemma \ref{lem:projequiv} we will show that $B^*$ is projective in $\Va^*$. Suppose that $A \in \Va$ and $p: A^* \onto B^*$ is an epimorphism. Let $\overline{p}: A \onto B$ be an extension of $p$ which exists by Lemma \ref{lem:extend*}. Since $B$ is projective, there is some $q: B \into A$ such that $\overline{p} \circ q = id_B$. Then, $\overline{p}\restrict A^* = p$, $q[B^*] \subseteq A^*$, $p \circ q\restrict B^* = id_{B^*}$, and thus $B^*$ is a retract of $A^*$.
\end{proof}
Collecting results we arrive at the main theorem of this section:
\begin{theorem}\label{thm: IA H pres} 
If $\Va \leq \IA$ is a locally finite $^*$-variety, then, $t(\Va) = t(\gamma(\Va))$.
\end{theorem}
\begin{proof}
By Theorem \ref{thm:full}, $\Va$ and $\gamma(\Va)$ are categorically equivalent, and thus, $t(\Va) = t(\gamma(\Va))$ by Theorem \ref{lem:Aequiv}.
\end{proof}
\section{Summary and outlook}

We have shown that the functor $\rho^*\colon \Lambda(\HA) \to \Lambda(\IA)$ preserves unification type, if $\rho^*(\Va)$ is locally finite. In the other direction, we have proved that $\gamma\colon\Lambda(\IA) \to \Lambda(\HA)$ preserves unification type, if $\Va$ is a locally finite $^*$ variety. The conditions we establish are sufficient for the preservation of unification type, but not necessary.

In future work we will investigate varieties which are not $^*$ varieties, in particular, we will consider locally finite varieties $\mathbf{W}$ of interior algebras of the form $\rho(\Va)$ for some $\Va \leq \HA$ which are in some sense a counterpart to $\rho^*(\Va)$, see the remark after Theorem \ref{thm:rho}. 

 \section{References}
\renewcommand*{\refname}{}

\end{document}